# Measuring Exploration: Review and Systematic Evaluation of Modelling to Generate Alternatives Methods in Macro-Energy Systems Planning Models

Michael Lau[1]

Neha Patankar[2]

Jesse D. Jenkins[1,3]

1. Department of Mechanical and Aerospace Engineering, Princeton University
2. Department of Systems Science and Industrial Engineering, Binghamton University
3. Andlinger Center for Energy and the Environment, Princeton University

## Abstract

As decarbonization agendas mature, macro-energy systems modelling studies have increasingly focused on enhanced decision support methods that move beyond least-cost modelling to improve consideration of additional objectives and tradeoffs. One candidate is Modeling to Generate Alternatives (MGA), which systematically explores new objectives without explicit stakeholder elicitation. Previous literature lacks both a comprehensive review of MGA vector selection methods in large-scale energy system models and comparative testing of their relative efficacies in this setting. To fill this gap, this paper provides a comprehensive review of the MGA literature, identifying at least seven MGA vector selection methodologies and carrying out a systematic evaluation of four: Hop-Skip-Jump, Random Vector, Variable Min/Max, and Modelling All Alternatives. We examine each method's runtime, parallelizability, new solution discovery efficiency, and spatial exploration in lower dimensional ($N \leq 100$) spaces, as well as spatial exploration in a three-zone, 8760-hour capacity expansion model case. To measure convex hull volume expansion, this paper formalizes a computationally tractable volume estimation algorithm using the sum of all combinatoric 2-d projections of higher dimensional convex hulls. Through these tests, we find Random Vector provides the broadest exploration of the near-optimal feasible region and Variable Min/Max provides the most extreme results, while the two tie on computational speed. We thus propose a new Hybrid vector selection approach combining the two methods to take advantage of the strengths of each. Additional analysis is provided on MGA variable selection, in which we demonstrate MGA problems formulated over generation variables fail to retain cost-optimal dispatch and are thus not reflective of real operations of equivalent hypothetical capacity choices. As such, we recommend future studies utilize a parallelized combined vector approach over the set of capacity variables for best results in computational speed and spatial exploration while retaining optimal dispatch.

# 1.0 Introduction

As decarbonization goals necessitate rapid, large-scale transformations of global, national, and subnational energy systems, macro-energy systems models are valuable tools to provide insights into the likely effects of novel technologies, investments, and policy interventions. Macro-energy system planning models are typically formulated as large-scale linear or mixed-integer linear programming problems with the explicit objective of minimizing cost, subject to various engineering, economic and policy related constraints. And yet, real-world decisions are almost never made entirely on a least-cost basis. Affordability is a central concern to most stakeholders, but far from the only salient objective or determinant of decision making. More nuanced objectives such as political feasibility, popularity, or energy security concerns frequently inform decision-making, while distributional outcomes shape stakeholder engagement in decision-making processes. Yet these multiple and often-competing objectives are typically omitted from explicit consideration in macro-energy systems planning models or treated only through limited scenario or sensitivity analyses.

Uncertainty in energy system models falls into two general categories. Parametric uncertainty describes uncertainty in the range of values for a given parameter or set of parameters. Due to the forward-looking nature of energy system models, many parameters—the cost of availability of various generation or storage technologies, changes in demand, or the price of fuels—will exhibit substantial parametric uncertainties. Structural uncertainties involve uncertainty in the construction of the model itself. Models cannot fully represent the real world, and thus inherently make simplifying assumptions. While expert opinion, model validation, and inter-model comparisons can strengthen confidence surrounding the accuracy of modeled scenarios, structural uncertainty will always remain as to the true level of representativeness present in any given model. In other words, parametric uncertainty relates to uncertainty about the state of the world we attempt to model, while structural uncertainty relates to uncertainty about how to represent that state of the world in an abstract model.

Modeling to Generate Alternatives (MGA), a family of multi-objective optimization techniques originally introduced by Brill et al. (1979) and first applied to the energy systems modeling by DeCarolis et al. (2011), offers a potential solution to two key limitations of current energy system models: structural uncertainty and non-modeled objectives/constraints. MGA assumes that the option space of real-world decisions is contained by the technical requirements of the system and a budget, defined as the maximum acceptable system cost. The methodology then aims to map the bounds of that constrained feasible space as thoroughly as possible, with the aim of giving decision-makers a sense of the real options available to them and enabling them to subselect between options based upon their stated and revealed preferences. In doing so, it also enables stakeholders to mitigate some forms of structural uncertainty, particularly in the interior direction. If, for instance, a constraint was omitted that cut out part of the feasible space within budget, it could be imposed on the mapped set of solutions without needing to re-solve the larger model. Given its potential utility, various algorithmic approaches to MGA have recently been applied to a variety of energy systems planning problems. Table 1 presents a listing of 17 papers that apply MGA in energy systems contexts and details the authorship of each paper, the year each was published, along with the specific method/s they utilize and a brief description of their respective applications.

As Table 1 illustrates, there has been a proliferation of different MGA algorithms applied to macro-energy system models. Yet prior to this work, there has been only minimal systematic comparative testing of the performance of major MGA methods. This paper fills that gap by conducting systematic testing of four prevalent MGA methods, documenting the results and providing usable recommendations for best practices when applying MGA methodologies to macro-energy system model applications.

We identified two prior comparative explorations of MGA method efficacy: Minkowski et al.'s paper on MGA algorithms in land use (2000), and Lombardi et al.'s paper on computational trade-offs in MGA (2022). Both leave significant gaps in the literature addressed by this work. The Minkowski et al. study, is applied to a small model with no exploration of dimensional scaling and each method is limited to finding only 15 to 22 solutions, which is unrealistically small in modern applications and therefore fails to capture the strengths of some methods tested. Lombardi et al. do apply a series of MGA methods to a contemporary macro-energy system model and generate 210 alternatives for each method, a reasonable number. However, Lombardi et al. do not test a variety of modern methods, particularly MAA and a parallelized Random Vector implementation, but rather a limited set of other weighting methods they constructed.[1] Notably, the methods tested by Lombardi et al. are sequential, requiring knowledge of previous solutions, which prevents parallelization. Furthermore, Lombardi et al. provide only minimal spatial exploration data, showing only a few metrics, and they do not analyze other key characteristics of each method including dimensional scaling and runtime. The dearth of testing of MGA methodologies in the macro-energy system context thus leaves significant uncertainty in the literature as to the performance of each method including the portion of the near-optimal feasible space each method captures in large models and the computational scaling performance of each method.

---

[1] This set includes HSJ, a proportional variant of HSJ they call relative deployment which adds the fraction of the max capacity for each resource in the previous solution to their respective weights, a randomized variant of HSJ that adds a random number between 0 and 100 to each resource's respective weight in the previous iteration, and a running average methodology that updates the weight based upon the average of previously found capacities for each capacity variable.

*Table 1: Energy System Model MGA Papers, Applied Methods, and Applications*

| Authors | Year | Method | Application |
|---|---|---|---|
| DeCarolis | 2011 | HSJ | Carbon Mitigation |
| DeCarolis, Babaee & Kanungo | 2016 | HSJ | Proof of Concept |
| Li & Trutnevyte | 2016 | Capacity Max | UK Electricity Sector Transition |
| Berntsen & Trutnevyte | 2017 | Random Vector | Energy Supply Scenario Diversity |
| Price & Keppo | 2017 | Distance | Uncertainty Exploration |
| Sasse & Trutnevyte | 2019 | Efficient Random Generation | Distributional Tradeoffs in Renewables |
| Nacken et al. | 2019 | Capacity Max/Min | Renewable Energy System Exploration |
| Jing et al. | 2019 | HSJ with Epsilon Constraints | Urban Energy Systems |
| Neumann & Brown | 2020 | Capacity Max/Min | Energy Supply Scenario Diversity |
| Lombardi et al. | 2020 | Relative Deployment | Renewables Siting Policy |
| Pedersen et al. | 2021 | MAA | Renewables Siting, Land Use, Equity, Emissions |
| Lombardi, Pickering & Pfenninger | 2022 | Relative Deployment, HSJ, Additional Random, Evolving-Average | Energy System Capacity Exploration and Methodological Comparison |
| Patankar et al. | 2022 | Random Vector | Land Use and Capacity Decisions in American West |
| Neumann & Brown | 2023 | Capacity Max/Min | Investment Decisions Given Price Uncertainty |
| Sasse & Trutnevyte | 2023 | Capacity Max/Min | Regional Interdependency in Europe's Electricity System Transition |
| Grochowicz et al. | 2023 | MAA w/ Chebyshev Ball | Finding capacity layouts for robust solutions to multiple weather year data |
| Pedersen et al. | 2023 | MAA | 55% Decarbonization Scenarios for European Electricity Sector |

This paper addresses this gap in the literature by reviewing the procedure and performance of four families of MGA algorithms: (1) the Hop-Skip-Jump method first proposed by Chang, Brill, & Hopkins in 1983 and applied to macro-energy systems models by DeCarolis (2011); (2) Capacity Max/Min first used by Li & Trutnevyte (2016); (3) Modeling All Alternatives (MAA) first described by Pedersen et al. (2021); and (4) Random Vector methods created by Berntsen and Trutnevyte (2017). Other sequential methods, including the Relative Deployment, Additional Random, and Evolving Average are not examined here as they can largely be considered subsets of HSJ, just with slightly differing weights (Lombardi, Pickering & Pfenninger, 2022). Efficient Random Generation, meanwhile, is a specific class of Random Vectors, but only specifies random weights for a subset of capacity variables in each run, leaving the others weighted at 0 (Sasse & Trutnevyte, 2023). Thus, it should be relatively well represented by the Random Vector implementation included here. Finally, the Distance algorithm is not included here due to the combination of its purely sequential nature and high computational intensity, as it requires the computation of Manhattan distance between each pair of points in the polyhedron. Furthermore, distance type methods have only been used in one paper that we found (Price & Keppo, 2017).

The literature review and testing of each method are carried out to answer the following key research questions about each of the MGA methods, namely:

1. How does this method explore the near-optimal feasible space?
2. How does this method perform from a runtime perspective?
3. How efficiently does this method find new solutions?
4. How much of the space does this method capture?
5. How does this method perform temporally and spatially as problem size increases?
6. How well do small-scale timing and spatial results apply to energy system models?

The rest of this paper is structured as follows: In Section 2, we provide a thorough description of each selected MGA methodology: Hop-Skip-Jump, Capacity Min/Max, Modelling All Alternatives, and Random Vectors. Section 3 explains and describes our experimental procedure, including devised metrics for outcomes of interest, testing protocols, and the testing setups used. Section 4 details the results of the testing, comparing the attributes of each method. Finally, Section 5 provides a discussion of the implications of those results for future studies using MGA, including best practices for implementation.

## 2.0 Algorithmic Description

Modeling to generate alternatives (MGA) encompasses a variety of methods designed to explore the near-optimal feasible space of a constrained optimization model. In the case of least-cost objective functions, MGA maps the feasible set of outcomes with tolerable costs (as specified by a budget constraint), which allows modelers and stakeholders to subselect modeled outcomes that satisfy a range of alternative objectives, including non-modeled but quantifiable outcomes, while emphasizing the degree of flexibility available in solution topographies. Further, in the context of energy transitions, proving the existence of many near least-cost alternative solutions can bring confidence that the path to decarbonization is not narrow. Thus, MGA methods are a powerful tool to address key weaknesses present in otherwise cutting-edge least-cost optimization models.

MGA methods consist of three primary steps, which may be accomplished in any number of ways. The steps are:

1. Solve the original optimization problem with a primary objective (e.g., minimize cost);
2. Transform the current objective into a new slack constraint added to the original feasible region;
3. Generate a new objective and optimize to identify a new, near-optimal feasible solution within the budget slack constraint;
4. Repeat Step 3 many times to generate a range of candidate solutions.

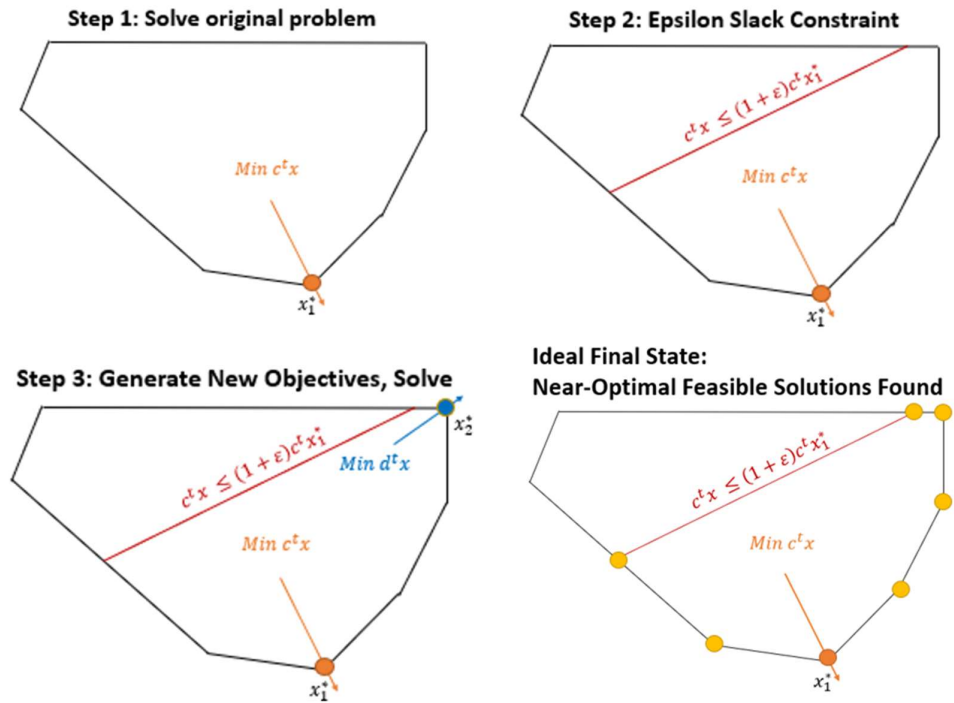

Figure 1: Stepwise Illustration of the Modelling to Generate Alternatives Procedure

These steps may be seen visually in Figure 1. The logical flow of MGA methods can be thought of as a procedural exploration of an unknown n-dimensional space, analogous to measuring the dimensions of a room in three-dimensional space. Every MGA method requires first finding an optimal solution, akin to stepping into the room we wish to measure and finding a corner to serve as a reference. Next, a target or slack constraint is imposed on the previous objective statement. In a least-cost optimization model, such a constraint would require all future iterations to fulfill all constraints in the model, as well as maintain an objective value within $x*(1+\varepsilon)$, where x is the cost optimal value and $\varepsilon$ a specified budget slack. The imposition of the slack constraint forces further exploration of the feasible space to remain close to the optimum, while allowing sufficient space to demonstrate all types of flexibility inherent to the given model structure. Returning to our room analogy, imposing the slack constraint is akin to adding a new wall that partitions the room, then exploring the corners of the newly constrained room. Finally, we measure the room.

Each method attempts to find the most diverse set of extreme points of the near-optimal feasible space possible, ideally employing an efficient and rapid procedure, in order to gain the fullest possible understanding of the extent of the near-optimal feasible space. For convex optimization problems, this is done by repeatedly creating a new convex objective function and solving the resulting problem, thereby finding a new extreme point of the convex hull of the optimization problem and expanding the known region of near-optimal feasible space. Holding computational power and choice of convex optimization solver constant, the way these new objective statements are chosen largely determines the performance of each sub-type of MGA procedure. Ideally, the algorithm is considered finished, and the near-optimal feasible space thoroughly explored, when the change in volume of each additional iteration converges, although in practice, many methods simply apply a pre-set number of iterations before terminating. Most existing algorithms subsequently prescribe some form of sampling procedure to represent the diversity of feasible solutions available.

MGA distinguishes itself from other multi-objective optimization approaches in two primary ways. Conceptually speaking, pure multi-objective optimization aims to find the Pareto frontier of a linear or mixed integer program with a set of different objectives. The Pareto frontier is defined as the set of feasible solutions of the model within which no objective function value can become more optimal (improve) without reducing another objective function value (Gunantara, 2018). The Pareto frontier represents the edge of the solution set in the objective space created by the set of objectives. MGA, on the other hand, aims to explore the near-least cost feasible space of the variable space which will include portions of the Pareto frontier with other objectives but may not fully represent it (e.g. it will exclude areas of the Pareto frontier outside the budget slack constraint). While MGA can be used to find an approximation of the Pareto frontier between many sets of objectives through ex post exploration of the set of MGA solutions, it does not explicitly search for these Pareto optimal solutions and thus may not include the optimal value for individual objectives within the objective set. The two can be combined to some extent through the addition of bracketing objective runs to any MGA protocol, which find the optimal value of secondary objectives in the near-least cost feasible space, similarly to the epsilon constraint method of multi-objective optimization (Mavrotas, 2009). See Figure 2 for a visual explanation of the differences between Multi-Objective Optimization and Modelling to Generate Alternatives.

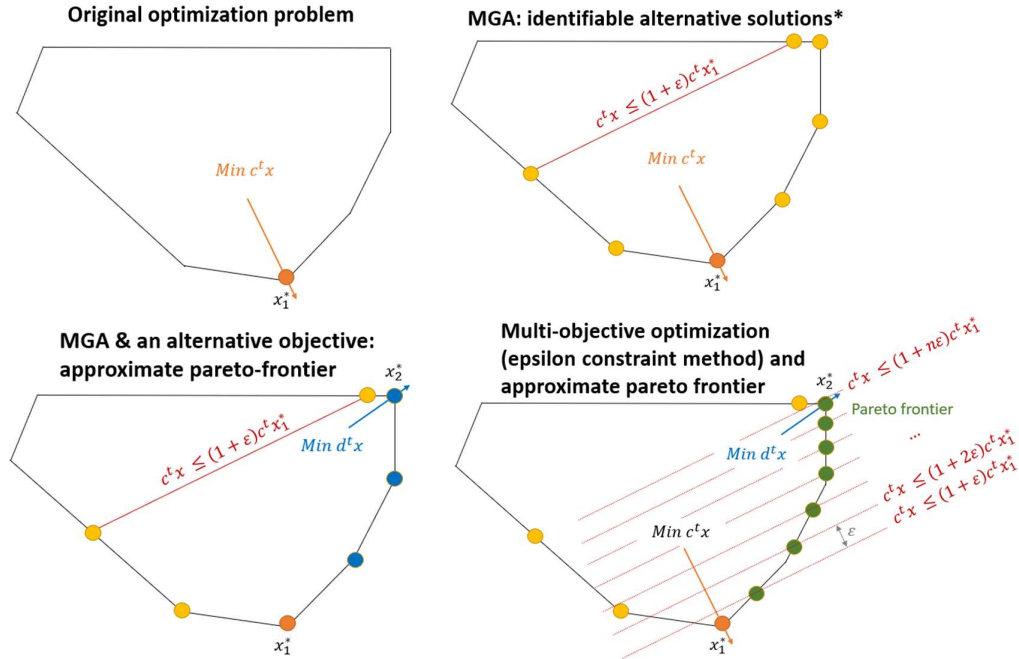

Figure 2: Illustration of the Differences Between MGA and MOO.

MGA can also help ensure the robustness of modeled findings to parametric and structural uncertainty. Figure 3 illustrates how various forms of parametric and structural uncertainty change the feasible space. Uncertainty in constraint formulation can take two forms – either the original set of constraints is overly constraining and a portion of the real feasible space is not captured by the model, or the original set of constraints is under constrained and a portion of the feasible space captured by the model is not truly feasible (e.g. the modeled constraints are a valid relaxation of the true constraint set). In the overly constrained case, where a modeled constraint is too tight or should be omitted, all solutions found are feasible but a wider range of feasible near-optimal solutions exists beyond the convex hull identified. In the under constrained case, where a constraint is omitted entirely or a modeled constraint should be tightened, some MGA solutions will no longer be feasible. If a portion of the real feasible space remains within budget, the solutions within that space will still be valid, meaning the true optimal solution and multiple near-optimal feasible solutions can be still found within the convex hull identified by MGA. MGA is least suitable for dealing with parametric errors in coefficients of modeled constraints that affect the shape of the feasible space, such as errors in transmission cost, losses, or demand and variable generation time series data. In these cases, some MGA iterates may remain feasible, but there is no consistent way to guarantee feasibility without recalculating the feasible space. MGA is most suited for dealing with parametric uncertainties affecting only objective function coefficients, as MGA involves modeling alternative objective functions subject to the original constraints.

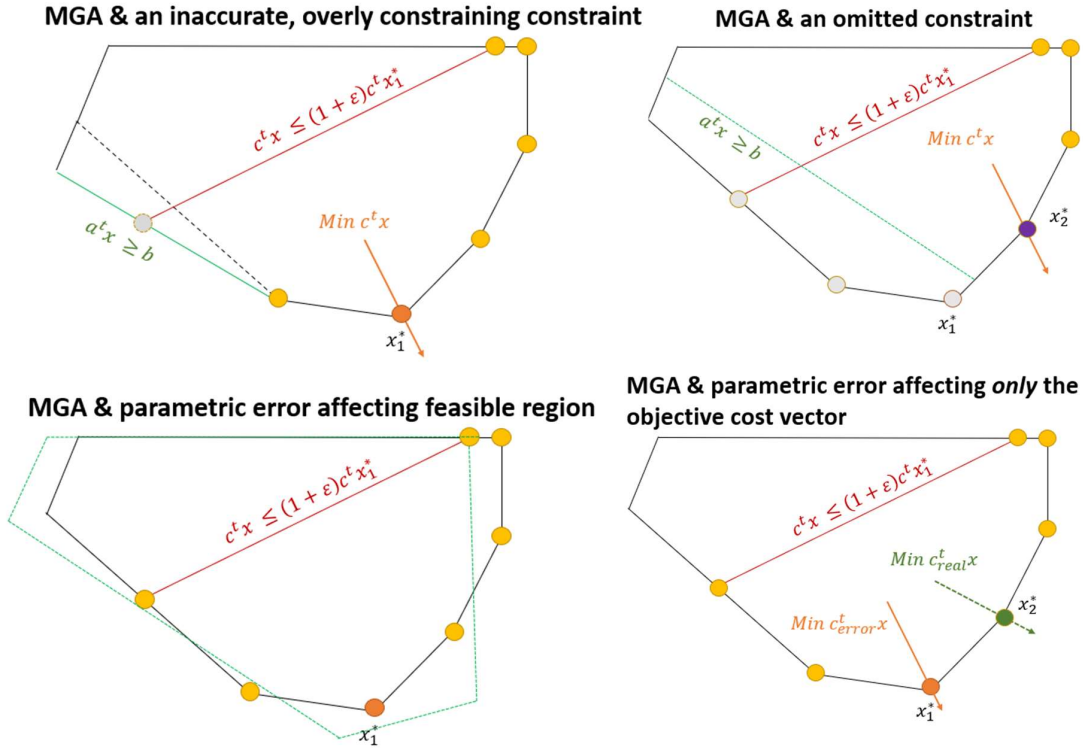

Figure 3: Illustration of the effect of MGA on various types of structural and parametric uncertainty

The following sub-sections describe each of the four specific MGA algorithms assessed in this study.

## 2.1 Hop-Skip-Jump MGA

Hop-Skip-Jump MGA (HSJ MGA) was formulated by Chang, Brill, & Hopkins (1983) and applied to energy system models by DeCarolis (2011). It consists of the following steps:

1. Find optimal solution to constrained linear optimization problem:
    a. $\min z = f(x)$
       Subject to: $A * x = b, x \in X$
2. Formulate the MGA problem:
    a. $\min p = \sum_{\{k \in K\}} x_k$
       Subject to: $A * x = b, f_j(x) \leq (1 + \epsilon) * z \; \forall j, \; x \in X$
3. Solve and repeat

HSJ MGA is distinguished from other MGA methods by its objective formulation. Within this formulation, $p$, the new objective, is the minimization of the weighted sum of all decision variables, with each weight determined by the number of solutions in which that variable has previously held a non-zero value (DeCarolis, 2011). As such, each iteration aims to explore a vertex of the feasible space which contains as few previously relied upon decision variables as possible.

While the HSJ MGA method has been shown to successfully generate a diverse set of scenarios (DeCarolis, 2011), HSJ MGA inherently limits its exploratory power due to its restrictive objective formulation rules. By only including the minimization of positive weighted sums in its objective function and iteratively relying on results of prior solutions, HSJ MGA tends to over explore the region of the space closest to the original objective, without deviating towards the rest of the polyhedron. HSJ MGA's strength is its focus on ensuring that it identifies edge cases that minimize each decision variable. By aiming to get only those solutions that expand edge cases for each optimization variable, it can begin to represent a fair breadth of the diversity of system configurations available while requiring relatively little raw computing power or sifting through a huge number of solutions, which can present challenges in cleaning and communicating findings. However, it will systematically under explore portions of the feasible region with higher values for several decision variables.

Limiting the search in this way inherently prevents HSJ MGA from comprehensively exploring the range of solutions available within the specified slack. The performance downsides of HSJ MGA are obvious even in a three-dimensional test system. To illustrate, the set of HSJ objectives and their optimal solutions within the 3-d test problem are written out in Table 2. HSJ identifies a total of four unique exterior points, out of eight total feasible exterior points in the full near-optimal feasible region for this 3-d test problem.

*Table 2: HSJ Example, all objective values negative to indicate minimization*

| Iteration | Objective Vector | | | Point | | |
|---|---|---|---|---|---|---|
| | x | y | z | x | y | z |
| 1 | -1 | 0 | 0 | 0 | 2 | 0 |
| 2 | -1 | -1 | 0 | 0.33 | 0 | 1.66 |
| 3 | -2 | -1 | -1 | 0 | 1 | 1 |
| 4 | -2 | -2 | -2 | 0 | 1 | 1 |
| 5 | -2 | -3 | -3 | 2 | 0 | 0 |
| 6 | -3 | -3 | -3 | 2 | 0 | 0 |
| 7 | -4 | -3 | -3 | 0 | 2 | 0 |
| 8 | -4 | -4 | -3 | 0.33 | 0 | 1.66 |
| 9 | -5 | -4 | -4 | 0 | 1 | 1 |
| 10 | -5 | -5 | -5 | 0 | 1 | 1 |

Given the scale of modern energy systems models, often with millions of decision variables and constraints, sequentially solving the HSJ algorithm is not a viable strategy to get a good approximation of the true range of system configurations present within the near-optimal feasible region. Authors have ranged in their estimation of what quantity constitutes sufficient solutions, but most have landed on a range of between 400 and 500,000 for a national-scale model. Time and computational constraints likely forbid most model applications from reaching even 400 solutions sequentially.

As HSJ MGA sequentially generates new objective functions, it is also not possible to parallelize the algorithm. As macro-energy systems modelling is frequently performed on super-computing clusters, parallelization of MGA algorithms may be achieved in one of two ways: by delegating individual model runs to distributed worker CPUs within a single job, or by superimposing the results of multiple jobs. When applying HSJ MGA, each objective function is formulated by checking which relevant decision variables have not yet appeared in a solution, then summing them, thus necessitating knowledge of the results of all previous iterations. The solving of each respective model therefore cannot be delegated to parallel processes within the same job. Furthermore, as HSJ is a totally deterministic process, creating multiple separate jobs would result in the same repeated output rather than additional value, making superimposing the results pointless. Distance-based MGA approaches, like those used by Price and Keppo (2017), have the same issue, reducing their efficacy in larger model applications.

HSJ has been successfully applied to the energy sector by DeCarolis et al. in 2016, along with a variant of HSJ using fractional weights determined by each generation technology's respective "contribution to electricity supply in the base case" (DeCarolis, 2016). Notably, due in part to the computational constraints of the method, only 16 MGA scenarios are examined, in addition to 3 traditionally derived scenarios (one base case and two CO2 cap cases). While HSJ is referenced by other papers in the energy modelling field, it does not appear in any other published studies found during the literature review here.

## 2.2 Modeling All Alternatives

Modeling All Alternatives, shortened here to MAA, is a methodology designed on the principle of outward expansion. Rather than seeking to get as far away from previous solutions as possible or incorporating new variables, MAA looks to maximize in new outward directions on every iteration. The method works by first finding a polyhedron of at least three initial exterior points, computing its convex hull, finding the face-normal vectors of the halfspaces that make up that convex hull, imposing each of them as the objective function of a cloned LP and maximizing the slack-constrained LP in those directions. The process is then repeated with the new set of exterior points. Expansion in this method happens in stages. In each expansion stage, the convex hull of the polyhedron formed by prior exterior solutions has an increasing number of faces, making it possible to parallelize an increasing number of cases with each stage, but it is necessary to collect the full set of exterior points identified by each parallel process before proceeding to the next stage, limiting the scale of parallelization possible. This algorithm proceeds with more expansion stages until either an iteration number is reached, or the volume of the convex hull converges. The steps are described programmatically below:

1. Find optimal solution to constrained linear optimization problem:
    1. $\min z = f(x)$
       Subject to: $A * x = b, x \in X$
2. Impose slack constraint
    1. $\min z = f(x)$
    2. Subject to: $A * x = b, x \in X$

    $f_j(x) \leq (1 + \epsilon) * z \ \forall j, \ x \in X$

3. Run another MGA method to find three or more exterior points to form initial polyhedron
4. Find halfspaces that make up the faces of the convex hull of the polyhedron through Qhull algorithm

5. Create new MAA problems and solve (ideally in parallel)
   1. $\min n * x$
   2. Subject to: $A * x = b, x \in X$
      $$f_j(x) \leq (1 + \epsilon) * z \ \forall j, \ x \in X$$
      Where *n* is the normal vector of a face of the convex hull (repeated vectors are discarded).
6. Collect solutions from 5 and repeat steps 4-5 until volume of convex hull converges or iteration count is reached.
7. Sample interior feasible near-optimal solutions through Delauney triangulation or Monte-Carlo Method (optional)

MAA has two major limitations. The first derives from its usage of the Qhull algorithm to calculate the convex hull and face-normal vectors of the polyhedron constructed in each expansion stage. In our testing, as problems increase in size past roughly 10 dimensions and 20 constraints and the numbers of solutions in the convex hull grows greater than 20, the Qhull algorithm slows down drastically, often failing entirely. This will occur at different points depending on the computational power of the computer involved, but at the scales typically required for capacity expansion models, the computational cost of the Qhull algorithm will fundamentally constrain the use of MAA. It is possible to reduce the dimensionality of capacity expansion models for MAA (or other MGA algorithms) by focusing on identifying variation in a limited number of decision variables, such as through clustering generation types across regions or only including certain metrics, such as generation capacity or overall transmission buildout in the MGA calculation. However, such simplification may lose much of the detailed value that MGA results offer in the process. A similar issue occurs when using the Delauney triangulation sampling method to identify interior point near-optimal feasible solutions, which requires complex calculations carried out by Qhull or another convex hull algorithm that does not work on larger dimensional convex hulls. As such, MAA is not suitable to explore a high-dimensionality near-optimal feasible region. It is worth noting that this is not a flaw fixable by simply changing vector calculation method while retaining the complex, geometrically motivated vector solution process of MAA and similar methods. In the testing run for this paper, we developed a similar method to MAA that did not use convex hull calculation packages but rather vector subtraction to generate generally face normal vectors. However, we did not see improved performance relative to MAA with this change and discarded the method.

The second limitation relates to all methods focused on maximizing in directions chosen by some variation of face-normal vectors. Methods that expand in this fashion work by generating vectors that point between previously found extrema, thus generating a greater density of vectors in already heavily searched directions. This strategy unnecessarily expends resources that would often be better served searching relatively unsearched sections of the polyhedron. While MAA attempts to minimize this difficulty by discarding vectors within a certain angle of one another, the sheer number of solutions present in even a 10-dimensional problem means that unique face-normal vectors could be generated almost ad-nauseum and find new solutions that do not violate the angular constraint but also do not contribute significant volume, expending valuable computational resources.

When MAA can compute all of the face normal vectors quickly, the consistent expansion in relatively separate directions allows it to find new exterior points on nearly every optimization. For small numbers of iterations on low-dimensional problems, this can be a computationally efficient strategy. As the dimensionality of a problem expands, MAA's two limitations kick in and the algorithm loses effectiveness, becoming completely unworkable for a certain dimensionality of problem. MAA has thus far been used in two papers, the one within which it was proposed and an exploration of the decarbonization of the European energy system (Pedersen et al., 2021; Pedersen et al., 2022).

### 2.3 Random Vector MGA

Where the preceding two methods base their new optimization directions on the results of previous iterations, Random Vector MGA takes a different approach focused on generating the most diverse set of optimization vectors possible, ideally evenly sampling all regions of the n-dimensional polyhedron through the generation of a large quantity of random vectors corresponding to the number of iterations desired. The algorithm is as follows:

1. Find optimal solution to constrained linear optimization problem:
    1. $\min z = f(x)$
       Subject to: $A * x = b, x \in X$
2. Impose slack constraint
    1. $\min z = f(x)$
    2. Subject to: $A * x = b, x \in X$

    $$f_j(x) \leq (1 + \epsilon) * z \ \forall j, \ x \in X$$

3. Introduce random objective vector and solve
    1. $\min g * x$
    2. Subject to: $A * x = b, x \in X$

    $$f_j(x) \leq (1 + \epsilon) * z \ \forall j, \ x \in X$$

    Where $g$ is a random n by 1 vector.

4. Repeat Step 3 a desired number of times or until the volume of the convex hull described by identified solutions converges.

In our testing, Random Vector MGA performs very well on problems with dimension greater than 10, on all metrics including runtime, unique solutions found, and unique solutions found per optimization. This is due to several factors. First, the only computationally intensive step of the algorithm is solving each LP; generating a new objective is trivial (e.g., involves only calculating a random vector), unlike MAA or other expansion-based algorithms which involve computationally expensive calculations to determine new search vectors. Random Vector MGA is also highly parallelizable due to its lack of memory. Separate Random Vector algorithm iterations can be implemented on all processors, limited only by the power and number of processors available.

The Random Vector methodology is not without downside. Since the objective of Random Vector can point anywhere in the MGA space, it rarely points down an axis, or in the space spanned by only a subset of axes. Thus, while it consistently performs well by all metrics, the Random Vector methodology can struggle to find the maximum and minimum possible feasible values for specific variables. The notable benefit associated with this problem is that each set of variables does not have to be intentionally maximized or minimized at some point in order to have flexibility in that dimension captured, since almost every vector will include some weight for every variable. As a result, the Random Vector methodology can be sure to capture flexibility in every dimension in a high dimensional problem in relatively few runs. We can search for more extreme values in a subset of variables by giving the others a weight of zero, but doing so removes any guarantee of capturing flexibility in every dimension. In the macro-energy systems context, this means that for problems with salient decision variables, such as a large number of choices for the location and type of generation or storage capacity additions or retirements, and when higher computational complexity for the original LP restricting the number of MGA iterates it is possible to compute, Random Vector will still be able to capture a reasonable estimate of the near-optimal feasible space for all dimensions.

## 2.4 Capacity Min/Max

Capacity Min/Max describes a family of MGA methods which rely on maximizing and minimizing a set of randomly or intentionally selected variables in each MGA iteration. Though similar to Random Vector MGA in concept, as they both maximize and minimize various combinations of variables, the Capacity Min/Max method differs from its counterpart in the way it formulates the weights assigned to each variable. Algorithmically speaking, when a set of MGA variables is presented, Capacity Min/Max assigns subsets of variables weights chosen from the set of integers [1,0,-1], indicating that that subset of variables should be minimized, maximized, or allowed to flex to accommodate other variables' optimization. Like the Random Vector methodology, Capacity Min/Max performs well on large problems due to its parallelizability and computationally inexpensive objective creation procedure. However, its exhibits the opposite characteristics as Random Vector with regards to dimensional inclusion and extreme points along axes: where Random Vector could explore all dimensions in relatively few iterations but fail to capture the most extreme points for specific variables, Capacity Min/Max tends to capture very extreme points for the variables included in each objective statement but will miss exploring dimensions that are not explicitly maximized or minimized in one or another objective statement. This problem can be blunted by ensuring all variables appear in at least one objective with a weight of 1 or -1. However, there is a trade-off here as well. When many variables are maximized or minimized in the same objective, they typically pull against one another, resulting in less extreme values being found for each individual variable. Thus, in high dimensional problems with high computational complexity and a smaller number of iterations possible in a given timeframe, attempting to maximize or minimize all variables can result in a similar shape as that from Random Vector. Balancing this trade-off is a key challenge of implementation for both the Capacity Min/Max and Random Vector methods, and we explore possible hybrid solutions combining these two methods in our experimental procedures.

## 3.0 Experimental Procedure

To achieve this end, we utilized three different testbeds: a three-dimensional known, visualizable model, an N-dimensional randomized scalable LP generator, and a test scenario from the GenX electricity system model (Jenkins et al., 2022). The three-dimensional model was designed and utilized to provide visual information as to the exploratory tendencies of each MGA method. By examining the objectives generated and points found by each method over many runs, we were able to gain greater insight into the way each method procedurally finds points, making their strengths and weaknesses more apparent than they are in higher-dimensional spaces where all points are not known, and the space is not completely visualizable. The N-dimensional testbed was designed to provide information as to how each method scaled with model size and increasing numbers of dimensions included in MGA objective vectors. While these models were no longer visualizable, the N-dimensional testbed allows for large-scale and repeated rapid testing of each method on sets of random models, generating the core of the temporal, spatial, and dimensional scaling data for this analysis. Finally, three of the methods, HSJ, Capacity Min/Max, and Random Vectors are implemented in a macro-energy systems modelling context, specifically a Three Zone ISO New England model implementation of GenX, a contemporary open-source electricity system capacity expansion model, to test their applicability to capacity expansion models.

### 3.1 Testing Regimes

To confirm that each method was implemented properly and was behaving as expected, a known three-dimensional problem was designed, shown in Equation 1.

$$minimize\ x_1 + 2x_2 + 2x_3$$

$$subject\ to: x_1 + x_2 + x_3 \geq 2$$
$$x_1 \leq 3$$
$$2x_2 + 3x_3 \leq 5$$
$$10 \geq x_1, x_2, x_3 \geq 0$$

Equation 1: Initial 3-Dimensional LP

As the three-dimensional space is small with relatively few vertices, we imposed a slack of 3, resulting in the shape visualized in Figure 1. We used the Qhull algorithm to find the full convex hull of the slack-constrained test problem for reference when verifying the implementation of each MGA method (Barber et al., 1996).

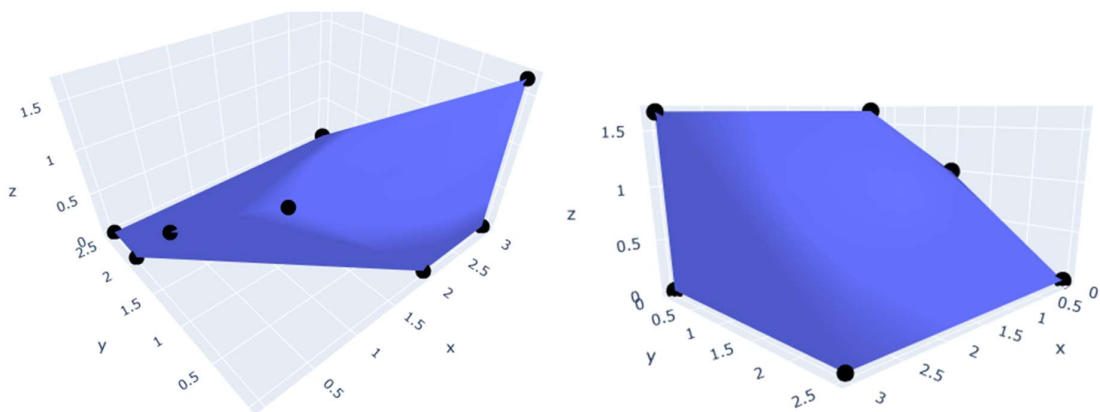

Figure 4: Three-Dimensional Test Problem Convex Hull

While the three-dimensional test provides visual feedback as to how each method explores a space, it cannot provide adequate information about each method's performance in larger, higher-dimensional spaces akin to those important to energy modelling. To address this limitation, we developed a set of testbeds. The test script first generates an LP with *n* variables and *2n* constraints. Each variable was additionally constrained to have a value no less than 0 and no greater than 10. The objective statement was randomly selected as were initial constraints.

$$minimize \sum_{j=1}^{m} c_j x_j$$

$$subject\ to \sum_{j=1}^{m} A_{ij} x_j \geq b_i \quad i = 1, \dots, 2m$$

$$10 \geq x_j \geq 0, \quad j = 1, \dots, m$$

Equation 2: Random dimensional LP formulation

We ran each MGA method on a series of random problems ranging from 3 to 1000 dimensions with a slack of 0.1 (10%). The full testing regime can be found in Table 2. In addition to evaluating runtime and parallelizability, MGA methodologies are evaluated based upon how well they explore the near-optimal feasible space, defined as the area of the convex hull created by the original problem and cost constraint.

To test this quality in an electricity system model, a number of trials were conducted within the GenX capacity expansion model, utilizing the three-zone ISONE example case that is part of the GenX repository (Jenkins et al., 2022). All computation for this project was done on CPU nodes on the Della computer cluster at Princeton University.

*Table 2: Testing Regime, all dimensions run 10x up until 100 dimensions, then 2x each*

| System | Method | MGA Dimensions | LP Variables |
|---|---|---|---|
| Randomized LP MGA Testbed | Hop-Skip-Jump | 3, 5, 10, 20, 50, 100, 1000 | 3, 5, 10, 20, 50, 100, 1000 |
| | Random (Sequential) | 3, 5, 10, 20, 50, 100, 1000 | 3, 5, 10, 20, 50, 100, 1000 |
| | Random (Multithreaded) | 3, 5, 10, 20, 50, 100, 1000 | 3, 5, 10, 20, 50, 100, 1000 |
| | Modeling All Alternatives | 3, 5, 10 | 3, 5, 10 |
| | Capacity Min/Max | 3, 5, 10, 20, 50, 100, 1000 | 3, 5, 10, 20, 50, 100, 1000 |
| GenX Example System | Hop-Skip-Jump | 18 | ~1,400,000 |
| | Capacity Min/Max | 6 | ~1,400,000 |
| | Random | 18 | ~1,400,000 |
| | Random + Capacity Min/Max | 18 | ~1,400,000 |

### 3.2 Convergence Metrics and Criteria

Stopping criteria and convergence metrics for MGA algorithms have not been discussed extensively in existing literature. Most papers utilizing the method tend to document the number of iterations completed, along with the number of points generated, without any guarantee that they have explored the majority of the space, instead relying upon ex-post methods to select a small subset of maximally diverse options from amongst the identified solutions. Pedersen et. al (2021), do attempt to find a convergence criterion for MGA methods, the only paper we identified in the macro-energy systems literature to do so. They attempt to use QHull to provide volume calculations, then measure the difference in value calculated each iteration, terminating the method when the difference in volume falls below a predetermined termination threshold. As there are pseudo-infinite vertices in high-dimensional convex hulls, like those created by macro-energy systems models, it is impossible to find every extreme point in reasonable runtimes. Thus, true convergence of the discovered convex hull with the true shape of the model is impossible. Conceptually, however, Pedersen et. al.'s implementation should give as close to a guarantee that the majority of the space has been captured as is possible. Unfortunately, due to the NP-hard computational complexity of evaluating high-dimensional convex hull volumes with many points, Pedersen et. al.'s implementation is not practical unless a small number of decision variables are included in the alternative MGA objective functions.

It is possible, however, to make similar convergence criteria work for complex models, with a procedure first used by Patankar & Jenkins in 2021 which we have termed Volume Estimation by Shadow Addition, or VESA. We achieve this by first generating all combinatorial pairs of MGA variables. We then project all solutions onto each of these subspaces and take the area of the resulting 2-D convex hull of each dimensional pair, then sum them all to represent the total volume. Using VESA, we can capture most facets of the relevant convex hull without much computational difficulty as each 2-d hull is essentially a slice of the high dimensional polytope in those two dimensions, similarly to a CT or PET scan.

VESA is not a perfect estimation of volume and will not capture points that fall within the "shadow" of the convex hull. An example of this case can be seen in Figure 5. As demonstrated, however, this shadow is limited in volume and orientation. Any point added outside of the blue shadowed area will be captured by our metric. Thus, most volume additions are captured. Since, by definition of convexity, convex hulls cannot shrink in volume with the addition of a point, VESA will catalogue most volume additions to the convex hull. The change in VESA over a set of iterations can then be measured and compared to a convergence threshold, allowing for a consistent and computationally tractable termination criteria which provides some reassurance that the method has found most of the volume that it will find, and that further iterations are not worth their computational cost (e.g. are no longer increasing the area of the explored polyhedron).

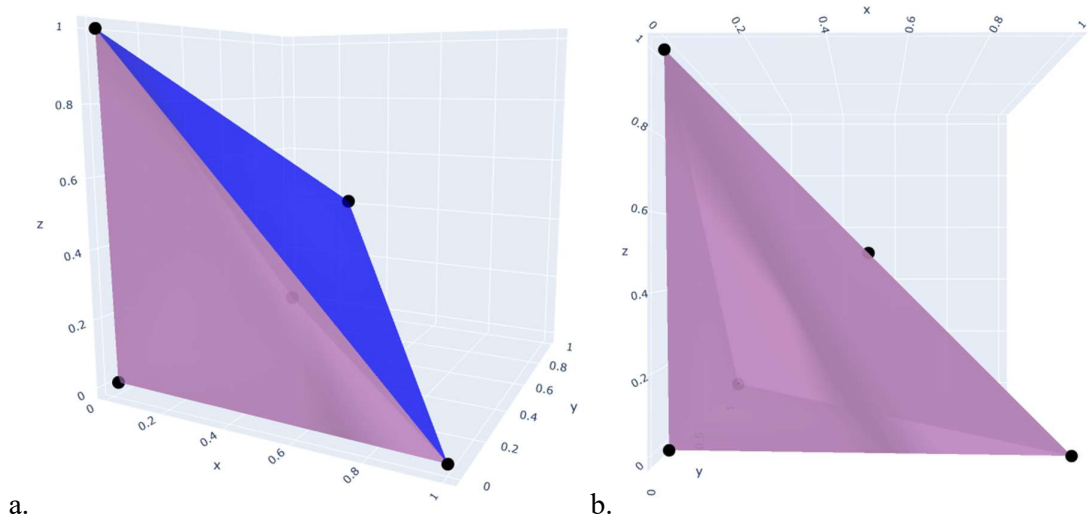

a.  b.

Figure 5: a. Three dimensional volume addition (blue) to tetrahedron (pink) and b. two dimensional slice with no area increase demonstrating that a point can add volume without adding to any crosssectional area

Two obstacles must be considered when using this termination criteria. First, when used as a convergence metric in multi-threaded implementations, different processors have only partial access to the set of previously identified solutions and thus cannot compute complete volume estimates, resulting in some variation in volume estimate between iterations depending on the solution sets available to each processor. This can be largely obviated by either implementing distributed arrays or running many iterations on each parallel processesor, returning results from this 'batch' and sorting the results to determine if convergence criteria have been met. If not, another set of parallel batches can be initiated. Second, methodological deficiencies may lead to the appearance of convergence, without fully exploring the space. This difficulty is most present with HSJ, which has very limited exploratory efficiency. It is not uncommon for HSJ iteration volume estimates to converge to relatively small values compared to the volume found by other methods in the same space, implying that HSJ is getting stuck in a region without a way to escape and diversify its vector selection. Adding a convergence metric will not change that and will only give a false sense that the algorithm has captured the majority of the space, when instead it has only found all of the solutions the method is capable of finding – two very distinct concepts.

### 4.0 Results

As discussed earlier, when considering the quality of an MGA algorithm, two primary characteristics are prized: the speed with which the algorithm can be completed relative to problem scale, and the exploratory power of the algorithm within that time (e.g., the ability to approximate the convex hull of the full near-optimal feasible region of the problem). The following results indicate that the multithreaded Random Vector and Capacity Min/Max algorithms consistently performs the best of all MGA methods for larger scale problems.

## 4.1 Computational Performance: Runtime

MGA subtypes are largely delineated based upon their respective objective vector creation methods with a wide variance in computational complexity between methods. To illustrate the differences, objective creation was timed for all methods using the Julia time_ns() implementation over the set of dimensions detailed in Table 1. Each dimension and method ran 10 times, with MAA stopping before the 20-dimensional problems (approaching the maximize feasible dimensionality for QHull computations). The 100-dimensional problems and 1000-dimensional problems were only run twice each due to runtime constraints. The results of this testing are presented in Figure 6.

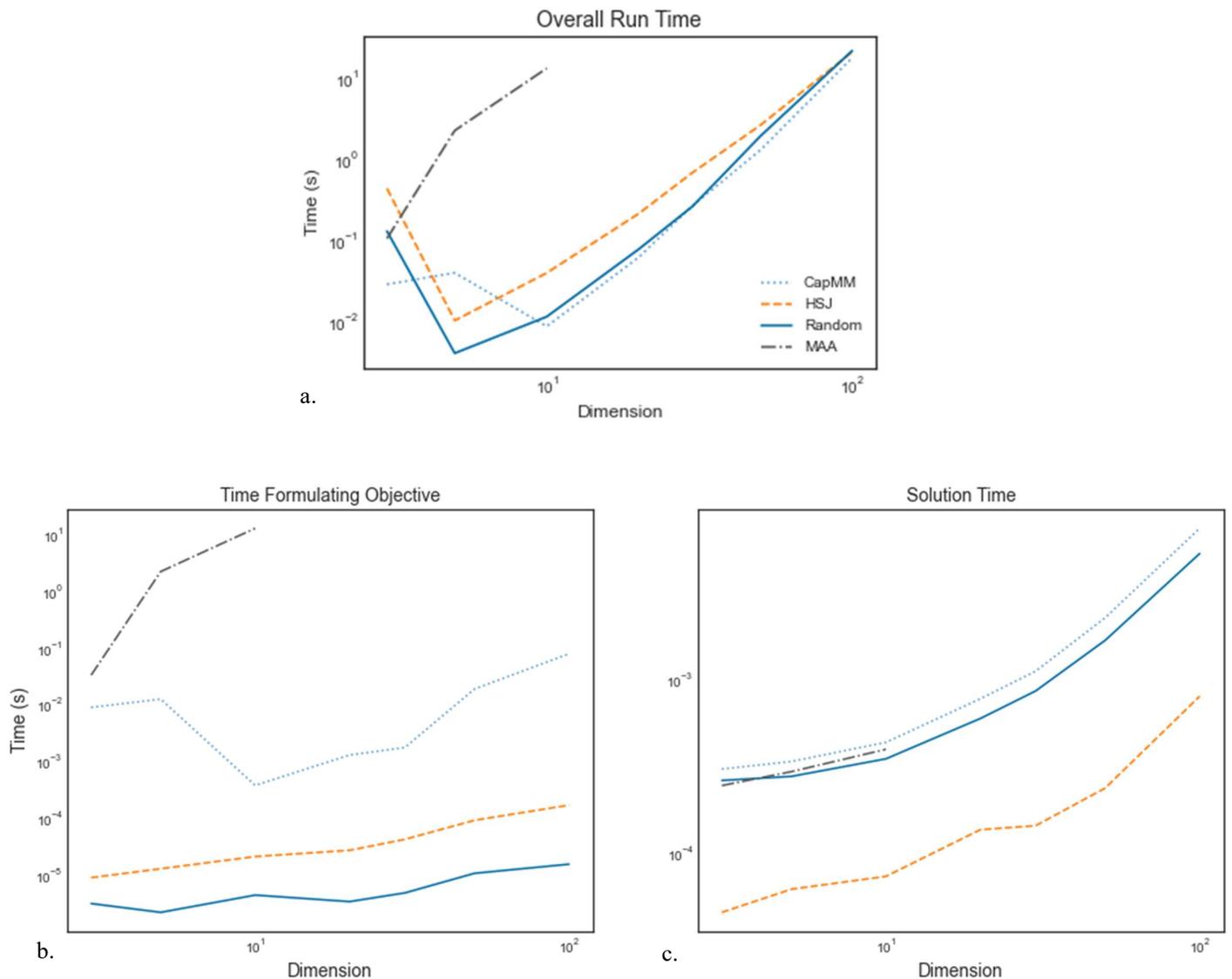

Figure 6: Key Runtime Metrics for all MGA methods
a. Average overall runtime for one run of all MGA Methods in randomized testbed, b. Max objective formulation time for all MGA methods, c. Average solution time for single iteration of all MGA methods in randomized testbed.

The results clearly divide MGA methods into two categories. Those with complex, geometrically based calculation for determining the next set of objective vectors, namely MAA, and those that rely on simple heuristics for the next set, HSJ, Capacity Min/Max, and the Random Vector method. Geometric methods tend to scale very poorly with dimension, due to the inability of any convex hull algorithm to calculate high dimensional convex hulls, which is confirmed in this case by MAA devoting most of its runtime calculating objective vectors rather than solving the LPs in question. Considering the dimensional scale of electricity system models, only electricity system models which abstract away from the zonal structure typically used, will be able to consider using MAA or similar methods. Including separate capacity decision variables for each zone within the MGA formulation, as seen in Lombardi et. al (2023), is out of the question for these methods as it dramatically increases the dimensionality of the MGA space explored. This is usually desirable in macro energy applications, as two solutions with identical total buildout of a given resource will be disparate in other metrics if the distribution of those capacity changes is different. Even in small-scale MGA with fewer than 10 explored dimensions, these two methods take significantly longer than their simpler counterparts, requiring significantly more calculation and allocated memory for any given number of iterations. Among the simpler methods, Capacity Min/Max took significantly more time to generate objective vector sets than HSJ and Random Vector. This difference is primarily due to one of two problems depending on problem size. In problems with many MGA dimensions, Capacity Min/Max requires a check to ensure that all variables of interest are maximized and minimized in one of the MGA iterations, otherwise those variables are frequently unexplored. In problems with few MGA dimensions, Capacity Min/Max requires checks for identical vectors since there are only three possible states for each weight, leading to higher likelihood of vector duplication. These additional checks take time and computational power relative to the HSJ and Random Vector methods. However, testing indicates that the computational time required for the Capacity Min-Max method to perform this check is relatively minor.

The simpler MGA methods generally performed well in terms of runtime, with HSJ rivaling the multithreaded Random Vector and Capacity Min/Max methods despite the latter group's ability to run eight models in parallel. Their convergence in times is largely due to the higher solve time the Random Vector and Capacity Min/Max methods experience with every iteration relative to HSJ. This difference can be attributed directly to the angular vector spread of the latter processes. By selecting points that lie on disparate sides of the polyhedron and rarely sampling similar directions sequentially, the Random Vector and Capacity Min/Max processes cannot readily warm start using previous iterations of the optimization solver as effectively as the less exploratory HSJ, resulting in longer solve times. When multithreaded, this issue is exacerbated, as each thread requires a separate instance of the model and solver, meaning that sequential warm starts cannot be carried over from run to run. However, as can be seen from the overall time graph, the speed gains associated with multithreading are enough to more than compensate for this issue and bring temporal performance comparable to HSJ to the Random Vector and Capacity Min/Max methods.

## 4.2 Spatial Performance: Unique Solution Set Size and Quality

MGA algorithm spatial performance is largely determined by two distinct metrics: (1) new solution efficiency, or the rate at which a new iteration of the model solution stage can be expected to find a unique solution; and (2) approximate convex hull volume, or an estimate of the extent to which the algorithm explores the near-optimal feasible region of the problem. Taken in combination, these two metrics provide an excellent means of determining the relative efficacy of each algorithm in discovering the important facets of the polyhedron in question as well as their relative efficiency in doing so. Both metrics are used to estimate the performance of each method in the N-dimensional randomized LP testing up to 1000 dimensions. It is worth noting, however, that the value of new solution efficiency decreases as the dimensionality and size of the space increases due to the exponential increase in the number of unique vertices present in the shape – several vectors separated by a hundredth of a radian may find unique solutions even though they may not contribute any significant additional volume to the convex hull. Accordingly, the MGA runs on the GenX test system are only evaluated by convex hull volume, rather than unique solutions produced.

### 4.2.1 Randomized LP Testbed

The spatial results of the Randomized LP testing program, detailed in Table 1, are reported in Figures 7, 8, and 9. Figure 7 relates both the total number of unique solutions each method was able to find in a set number of iterations and the efficiency with which each algorithm successfully found unique solutions (new solution efficiency). Figure 8 shows the portions of the 3-dimensional test convex hull discovered by each of the four methodologies discussed in this paper. Figure 9 illustrates the estimated convex hull volume found by each method for a given problem in 10, 100, and 1000 dimensions respectively.

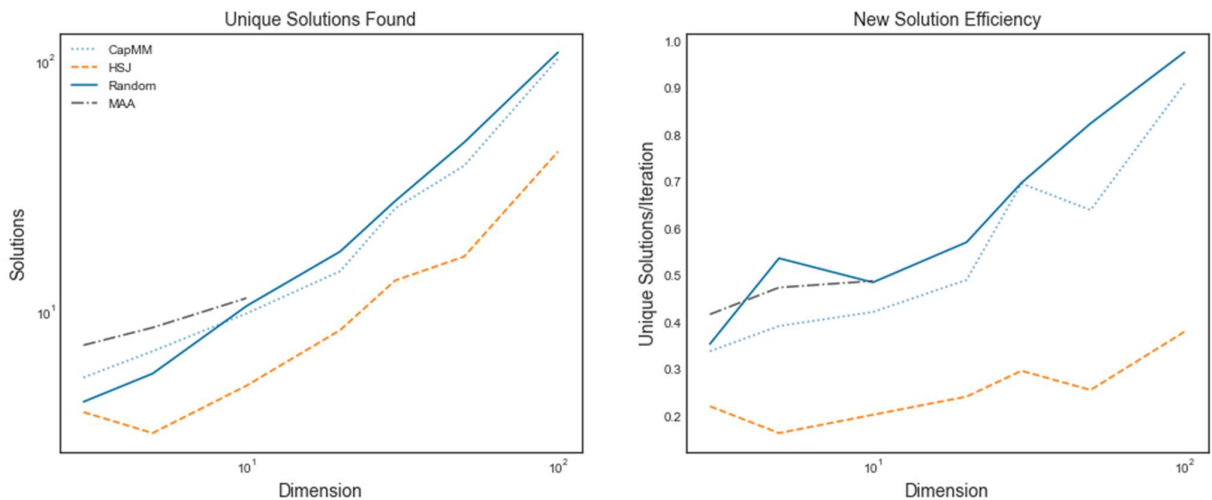

Figure 7: N-dimensional Randomized LP Unique Solution Performance, N = 10 for dimensions 5, 10, 20, 30, 50, 100, N = 5 for dimension 1000

Within smaller problems, where all methods were usable, MAA, Capacity Min/Max, and Random Vector consistently averaged more unique solutions identified and greater new solution efficiency than HSJ. MAA and Random Vector consistently had the highest new solution efficiency, likely due to the relatively precise geometric sampling method of finding face-normal vectors for MAA and the uniformity of sampling in Random Vector. Capacity Min/Max also performed relatively well in both efficiency and aggregate points found. For small problems, though, it was not quite as efficient as Random Vector and MAA. This is attributable to its tendency to search in lower variable subspaces, such as the x-y, y-z, or x-z plane in a 3 dimensional problem, which can risk missing points in which larger numbers of dimensions are non-zero, a problem that is less noticeable, though still existent, in larger dimensional problems.

As shown in Figure 8, MAA regularly finds all of the exterior points of the 3-dimensional test problem. There is a certain element of luck involved with MAA, however, as it relies upon a few initialization optimizations to create an initial set of points from which to generate its first round of face-normal vectors; these are obtained via the Random Vector method. Since MAA is generally restricted to smaller problems and the initialization step is meant to be short, relatively few Random Vector iterations are performed (usually around 2 iterations for the purposes of this study), and thus there is no guarantee that repeated explorations of the same space will reach the same spatial success or runtime every time. As the remainder of the MAA method is deterministic based upon the polyhedron described by the model in question and the quality of the starting points, the success of this initial step can have outsized effects on the quality of MAA search.

The Random Vector and Capacity Min/Max methods regularly perform the best of all methods on spatial metrics, retrieving the greatest number of unique solutions with the greatest efficiency for test problems with dimensions between 10 and 100, the range in which MAA has become computationally intractable and the scale of the space has not grown to the size where functionally any angular change results in a new solution. Heuristic-type methods like Random Vector and Capacity Min/Max are not as effective as MAA at finding new points in smaller spaces of dimension 10 or less due to the relatively large angular sweeps associated with each unique solution, which tend to reward geometrically derived vectors rather than the random or semi-random sampling associated with the heuristic-type methods. Despite not being as efficent as MAA in small environments, Random Vector and Capacity Min/Max do perform more than adequately in these smaller environments, consistently outperforming HSJ.

While the Random Vector methods does uniformly and efficiently explore the space, which works well for randomly selected polytopes, the engineered world tends to design solutions to be present along certain axes, and we tend to care more about solutions that reside along certain axes or metrics: i.e., lowest carbon emissions, lowest cost, greatest renewable energy share, or least transmission outcomes, etc.. As it is completely randomized, the Random Vector method does not typically capture these extremes; even in a 3-dimensional space, for instance, the random generation of a [1, 0, 0] search vector is relatively unlikely. This one notable weakness of the Random Vector method is clearly illustrated in Figure 6, where the Random Vector method fails to capture the (3, 0, 0) point that MAA successfully captures. Conversely, Capacity Min/Max is designed to explore these "down-axis" points extremely well, as it can only explore axes or lower-variable spaces, providing essentially the other side of the coin of Random Vector. This propensity is demonstrated in Figure 6 as well, as Capacity Min/Max specifically captures all points that are down an axis or directly between two axes, but misses points that do not optimize an axis or sum of axes, like (0, 2, 0), which is obscured by (0, 2.5, 0) to the [0, 1, 0] search vector, or (0.33, 0, 1.66). As will be discussed further in the Discussion section of this paper, these relative weaknesses of the Random Vector and Capacity Min/Max methods can be remedied by combining the two, resulting in a generalized Hybrid-type method combining the variable-including robustness of the Random Vector iterations, with a set of Capacity Min/Maxing

'bracketing runs' which serve to find extreme points by following objectives aligned with certain axes or outcomes of interest that may not otherwise be optimized by the Random Vector algorithm.

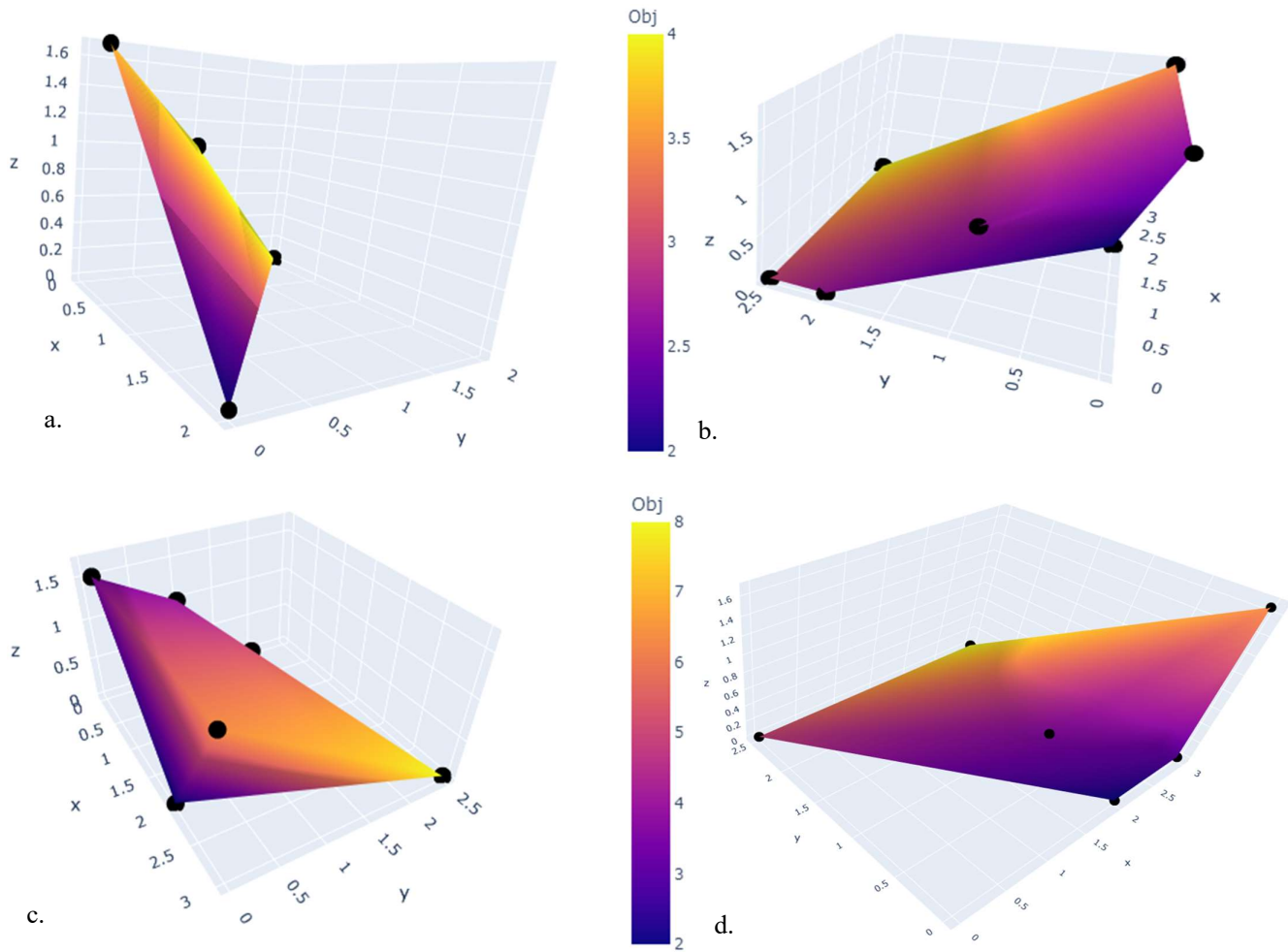

Figure 8: 3-dimensional example LP convex hulls found by each MGA method
a. Hop-Skip-Jump, b. Modeling All Alternatives, c. Randomized Vector, d. Capacity Min/Max

In addition to spatial accounting based upon the number of unique points found, and each MGA method's new solution efficiency, we consider spatial success for MGA methods in terms of their ability to identify solutions that describe a larger volume of the total near-optimal feasible region. Volume estimates for the high-dimensional polytopes found by each method were obtained through VESA, described earlier. While the volume estimate created by this method would not match the true volume of the n-dimensional convex hull, it should increase most of the time when the volume of convex hull increases thorough the inclusion of the new volume in one or more of the two-dimensional convex hulls. Thus, since we are applying this volume estimator consistently, any inaccuracies in volume estimation should be propogated across all methods, resulting in comparable measurements.

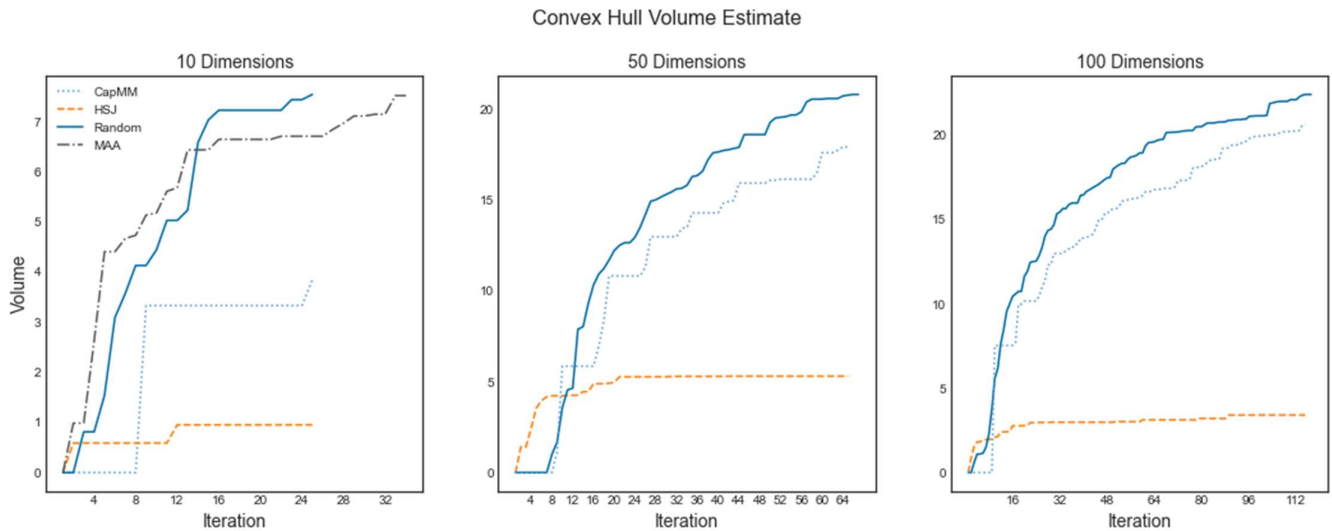

Figure 9: N-dimensional randomized LP convex hull volume estimates

In cases with 10 or fewer dimensions, MAA and Random Vector tend to converge to the same general volume estimation at roughly the same rate, indicating that they have discovered most of the space. As demonstrated in Figure 9a, for instance, the Random Vector and MAA methods converges to roughly 90% of the volume of the space in question in 12-25 iterations on average. The slow increase in volume as MAA begins to converge without finding significant new volume, from iterations 12 to 27, is due to the way MAA aims to split faces with objective vectors. This results in quick early expansion, but as the number of points increases, MAA generates many optimization vectors in the region with the highest density of previously found points, resulting in more inefficient iterations that capture little or no additional volume as demonstrated graphically in Figure 10. While these can be discarded in later iterations so they are not re-searched, they will occur in geometric algorithms of all types. Unlike MAA and Random Vector, Capacity Min/Max tends to struggle in smaller spaces due to its difficulty finding off-axis or off-plane points, resulting in substantial volume missed. However, it is still much better than HSJ, which does not perform well in volumetric terms in either these small dimensions or larger dimensions.

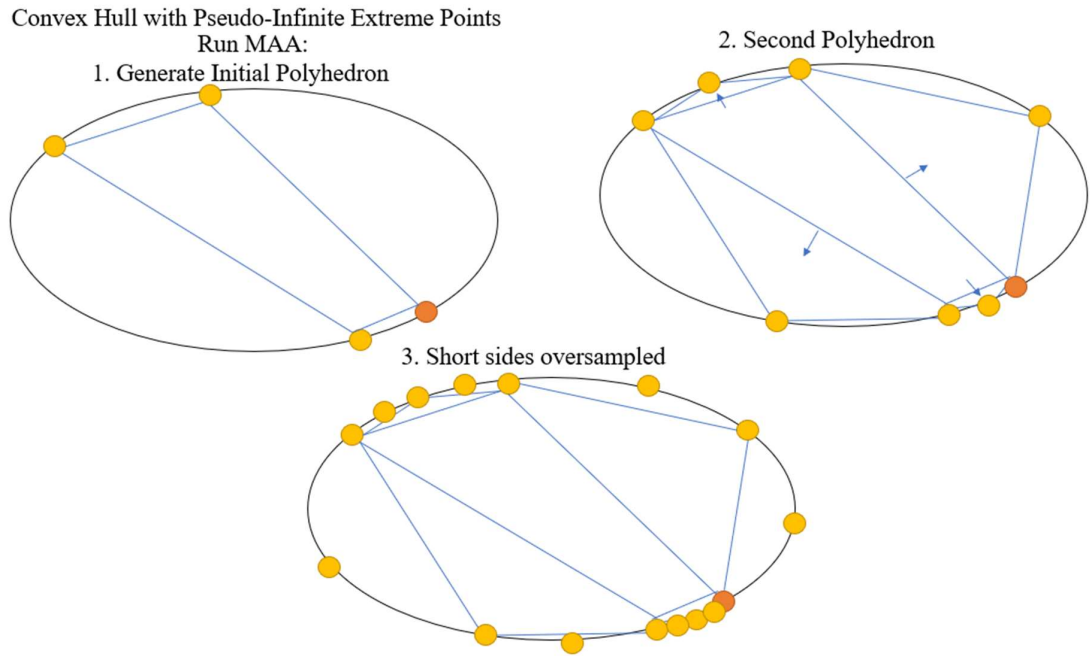

Figure 10: MAA tendency to oversample short sides of initial polyhedron

In the higher dimensional cases, MGA method options are limited to Random Vector, Capacity Min/Max, and HSJ due purely to computational constraints experienced by MAA. HSJ continues to struggle volumetrically in these larger problems, with its performance worsening especially in larger trials. In the 100-dimensional case trials, HSJ consistently failed to generate any meaningful volumes, converging around a total volume of 2 in over 100 iterations. By comparison, the volume found by the Random Vector and Capacity Min/Max algorithms typically require just 10 iterations or less to reach the maximum volume identified by HSJ, and these algorithms continued to increase the volume explored in subsequent iterations. Figure 9 gives an example of the difference in volume found between these methods in high-dimensional spaces. This is largely due to HSJ's inability to generate a broad diversity of objective vectors, tending to find solutions in the same general area of the polytope, as illustrated in both Table 1 and Figure 8. When the space in question expands to a very large number of dimensions, the number of vertices increases exponentially, reaching numbers where it is impossible to calculate all vertices. Thus, the number of unique solutions in any given region of the polytope similarly outstrips the number of optimization problems we can reasonably compute for a given question. Since HSJ does not explore in any geometrically-informed sense, but simply minimizes the variables which have already appeared, it tends to iterate around this same region, finding new points with different variable combinations that do not contribute much if any new volume.

Random Vector and Capacity Min/Max, on the other hand, tend to perform best in these larger problems. While it always performs relatively well regardless of dimensionality, higher dimensional problems play to two of the strengths of the Random Vector algorithm as they require more iterations and have a lower likelihood of coincident solutions. As the Random Vector algorithm functionally generates angularly uniformly distributed vectors, having a larger number of samples ensures a uniform distribution, while any angular perturbance will result in an entirely separate solution due to the sheer number of possible solutions present. In higher dimensional spaces, Capacity Min/Max typically generates a sufficiently large number of semi-random axis combination vectors that it also explores the space fairly well. While it still will not necessarily capture points between the lower-variable spaces it explores, the impossibility of capturing all points for any method makes this failing less apparent by comparison.

### 4.2.2 GenX ISO New England Test System

While the previous testing has demonstrated clear differences in performance between MGA methods when applied to a randomly selected LP with standard structure, we have not yet demonstrated that similar findings hold when applied to electricity system models. Regardless of their specific parameters, most linear capacity expansion models have a particular structure that is not at all spherical, as different classes of decision variables consistently vary in size by orders of magnitude. Total generation capacity and power flows between two nodes in each timestep are both variables present in almost any capacity expansion model, for instance, yet their numerical values will not be similar at all. When combined with relatively similar implementations across linear models for certain engineering requirements, like conservation of energy or estimation of transmission expansion cost, the overall structure of these models is often quite similar, likely resulting in relatively similarly shaped near-optimal feasible spaces, even if the specifics are different. This section aims to test the consistency of our findings for typical electricity system capacity expansion problems, namely that the Random Vector and Capacity Min/Max methods outperform HSJ on larger problems, as well as to demonstrate the theorized strengths of a hybrid method combining Random Vector and Min/Max searches. We performed a series of test runs with Random Vector, Capacity Min/Max, a combination heuristic-type method, hereafter referred to as Combo, and HSJ MGA in a well-documented Three Zone ISO New England test problem built into the repository of the open-source GenX electricity system model (Jenkins et. al, 2023). Note that as MAA and BA are computationally constrained to only very low dimensional MGA implementations, neither method was included in this set of tests.

The findings from this set of tests are presented in Figures 11 and 12. Figure 11 illustrates, in real capacity, cost, and emissions terms, the difference in near-least-cost feasible space captured by the Random Vector method and HSJ. Put another way, it illustrates the difference in the range of capacity options and cost and emissions outcomes captured for the test model running by each MGA methodology. Each subplot within the figure shows a specific pairwise tradeoff-space of feasible options, meant to demonstrate the level of flexibility and tradeoffs inherent between two resource types within a 10 percent budget slack. As the goal of MGA is to provide the best possible exploration of the near-optimal feasible space, the more area captured by a given method in each graph, the better it has performed, as more options and tradeoff relationships would be demonstrable to decisionmakers. Both consistency, as measured by finding volumes in all pairwise dimensions, and volume discovered are important in this context. Volume discovered is not sufficient if significant portions of the feasible space are missed in key dimensions, and exploration of all dimensions without significant volume is not sufficient either. Figure 12 collapses all these insights into one convex hull volume estimate, using VESA as discussed earlier.

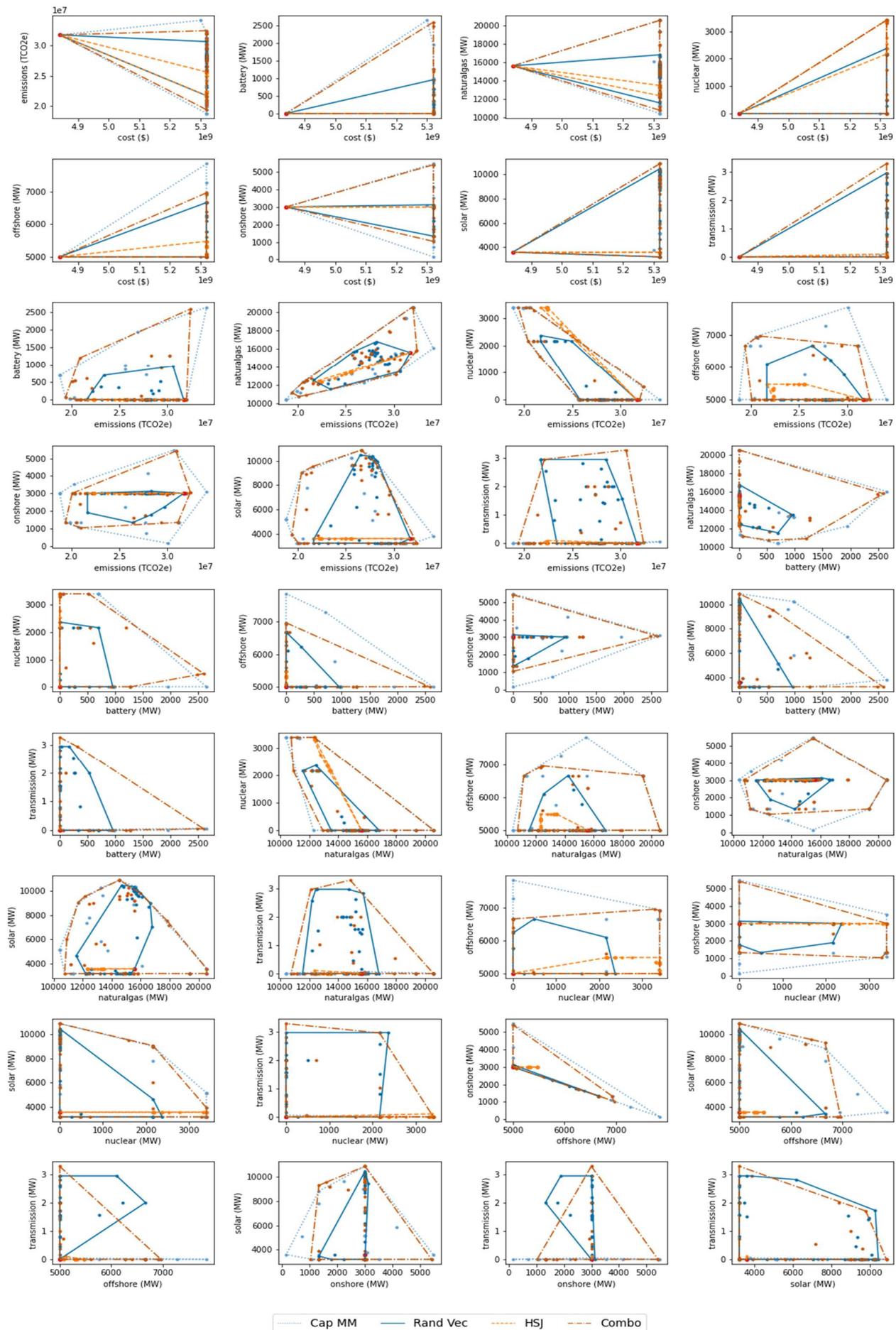

Figure 11: Comparison of pairwise convex hulls of GenX ISONE Three Zone test system captured by all methods, with all points. Optimal solution in red.

In all tradeoff-spaces visualized in Figure 11, the Random Vector method, Capacity Min/Max, and Combination significantly outperformed HSJ, capturing more volume and thus giving better vision on the level of flexibility present within the system, as well as the tradeoffs involved in decision-making around that flexibility. The additional volume captured has real ramifications for decision support. In nearly every case, the heuristic-type methods present a breadth of choices, showing flexibility in system composition and topology that would not have been known without such a high-quality exploration. Consider, for instance, the tradeoff plot between solar PV and offshore wind generation capacity. With HSJ, any decision-maker provided this modelling would conclude that while offshore wind has some capacity flexibility, as it can be built out to 5000 to 5500 MW in this specific system, while the system requires 3500 MW of solar PV capacity. However, with the improved visibility provided by, for instance, the Combo method, we can see that such a conclusion would be badly misinformed: there remains, in fact, substantial flexibility in both offshore wind capacity development and solar buildout within the 10 percent budget slack, with offshore wind capacities ranging from roughly 5000 to 7000 MW and solar capacities from about 3200 to 11000 MW, with combinations of capacities with high offshore wind and high solar or low offshore wind and low solar simultaneously achievable within the near-optimal feasible space.

The trends identified and theorized earlier among heuristic-type methods, namely Random Vector, Capacity Min/Max, and Combo, continue to evidence themselves in this applied testing. Random Vector shows great consistency across dimensions and trials, finding decent volume in all dimensions across all runs due to its inclusion of all MGA variables in each objective vector. Capacity Min/Max regularly finds the most extreme points for the variables included in its objective statements as it specifically mins and maxes them in limited combinations but can miss out on variables that are not specifically included. Combo shows the desired characteristics of both methods, with improved consistency for variables not specifically included in objective vectors relative to Capacity Min/Max and improved extremity and diversity of points identified relative to Random Vector. The difference in how each methodology treats technologies not specifically included in the MGA procedure can be visualized by examining the subplots that include transmission. Transmission capacity expansion was intentionally not included in the set of MGA variables here but is obviously an important outcome of interest. Notably, while Random Vector and Combo capture significant range in the transmission capacity dimension, Capacity Min/Max fails entirely. As the MGA dimension increases in applied studies, this will be true of more and more dimensions or outcomes of interest, with Capacity Min/Max thus missing important flexibility due to its hyper-focused nature. The Random Vector and Combo methodologies importantly mitigate this failing, albeit at the cost of some extremity of points explored.

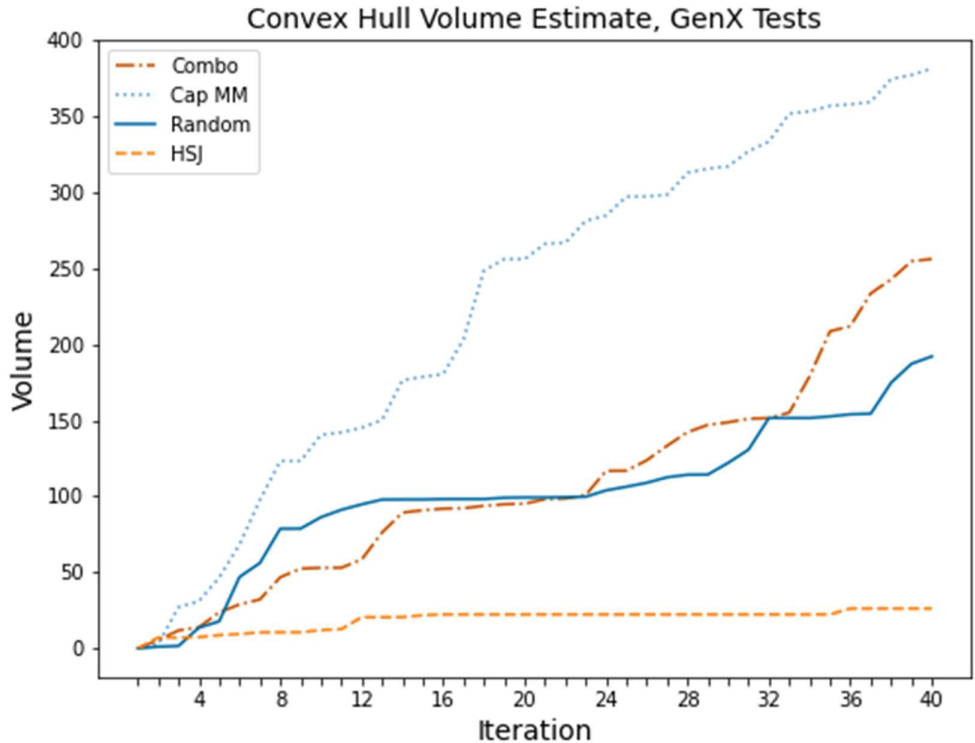

Figure 12: GenX ISONE Three Zone Test System Convex Hull Volume Estimate Comparison

We do see new behavior relative to the randomized LP testing emerge in the electricity system model context. Capacity Min/Max, which had underperformed relative to Random Vector in volume metrics in the randomized LP testing, captures significantly more volume than the Random Vector method, as shown in Figure 12. This is due to several factors. This trial run only includes 59 capacity decisions so the number of MGA dimensions present is relatively low, and 40 iterations were run due to the small size and relatively quick runtime of the model. Thus, Capacity Min/Max was able to specifically min and max nearly every pair of dimensions, resulting in high consistency relative to the Random Vector. Environments will become significantly less favorable to this style of method as the number of decision variables included in the MGA objective function increases, resulting in more MGA dimensions and longer runtimes. Random Vector still performs better than HSJ, but its low volume relative to Capacity Min/Max does highlight how the specific optimization of limited technology combos results in more extreme points. The Combo method, as noted earlier, does explore more volume than Random Vector (though less than Min/Max), while preserving its consistency in dimensional exploration.

## 5.0 Discussion and Conclusion

### 5.1 Implications on Suggested Usage

These findings carry several major implications for future MGA analyses, particularly in the macro-energy systems modeling domain, including insights related to both solution method and the interpretation of results from various methods used in the past. First and foremost, the Random Vector and Capacity Min/Max methods are evidently superior to other MGA methods from both a runtime and spatial perspective, and they can be easily combined to access both of their relative benefits (as in the Combo method presented in Section 4.2.2). Between their computational ease, simplicity of implementation, and parallelizability, these methods offer the lowest additional computational and coding burden of any of the MGA methods examined here. Random Vector and Capacity Min/Max or a combination of the two should work similarly well for any linear capacity expansion model and could easily be extended to mixed integer problems with no modification, albeit while incurring significantly longer solution time per iteration required to solve a MILP. Additionally, due to the low computational difficulty associated with calculating its objective vectors, the dimensionality of the MGA space can increase almost arbitrarily, the main drawback being that it will require significantly more iterations to gain any sort of volumetric conversion due to the sheer number of vertices in extremely high dimensional spaces.

      The strong performance of the Random Vector, Min/Max and Combo methods explored herein can unlock new types of analysis, including reducing variable clustering. However, expanding the types of variables included in the MGA objective function (e.g. the search space) beyond clustered regional capacity to annual generation outcomes or other generation-based metrics is not advised. Our testing has revealed that while least-cost dispatch is enforced by optimizing MGA objective functions over capacity variables, the same is not true of generation variables. When generation variables are explicitly included in the MGA objective function, these variables will be assigned alternative coefficients unrelated to the actual cost of generation. As a result, the optimizer tends to produce dispatch results that are not consistent with least-cost economic dispatch principles, resulting in unrealistic outcomes and unrealistic expectations for the performance of the proposed system in a real operations. This finding affects all metrics which are dependent on operations variables including but not limited to carbon emissions, utilization rates for thermal plants, operational costs, and local air pollution. An example of this issue is demonstrated in Figure 13. To demonstrate this phenomenon, we ran the same 3-zone example system with a random MGA objective on the annual generation for each type of generation and a budget constraint set up, then fixed the capacities found and reoptimized dispatch while minimizing cost of operations. To our knowledge, this shortcoming of traditional MGA approaches has not been noted in previous literature. Examples of papers using generation variables in their MGA formulation includes but is not limited to Berntsen & Trutnevyte (2017), and DeCarolis et al. (2012).

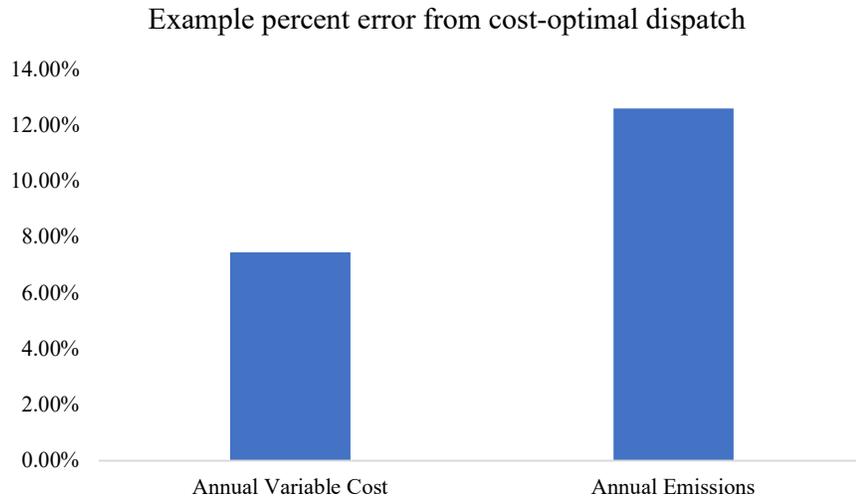

Figure 13: Example percent error from cost-optimal dispatch with fixed capacities in annual variable cost and annual emissions for MGA run on annual generation variables

It is also important to note that even though Random Vector and Capacity Min/Max performed extremely well in all tests carried out for this study, they can be parallelized further with the delegation of more computational resources. Increasing the quantity of iterates calculated will further improve the solution sets discovered by each of these two methods. As Random Vector is random, any two parallel set of runs can be superimposed to discover new solutions, a quality not associated with more deterministic methods like HSJ. Thus, in the same wall-clock runtime, the number of solutions that Random Vector could discover is theoretically limited only by access to computational resources. Especially in smaller dimensional problems, the same is not necessarily true of Capacity Min/Max. Due to the limited number of possible states for each MGA variable, the combinatoric likelihood that an objective vector is repeated increases as the number of iterates increases and the number of MGA variables decreases, making it important to check for duplicate vectors when superimposing multiple sets of separately generated Capacity Min/Max objectives or solutions to avoid wasting computational resources.

As mentioned throughout this paper, Random Vector and Capacity Min/Max can further be improved through their combination: adding lower-variable space specific runs to a set of Random Vector objectives either with the idea of 'bracketing runs', a preset list of objective functions determined to be of interest to the research question being explored including, for example, maximum solar capacity, minimum transmission capacity, or lowest carbon emissions, or through the Capacity Min/Max process of generating unique vectors with randomized weights selected from the set [-1,0,1]. These sorts of combination methodologies allow for the broad-spectrum exploration of the near-cost-optimal feasible space through the Random Vector method while retaining some ability to find extreme points in directions or lower-variable subspaces of interest. We demonstrate the strengths of this combination in Section 4.2.2, where the Combo method equals Random Vector in consistency and outperforms it in volume, and outperforms Capacity Min/Max in consistency while capturing slightly less volume. Combination methods like this appear to be the most effective at managing the tradeoff between extremity and consistency, and we thus recommend these as the best MGA methods currently available for use in macro-energy systems contexts (or similar domains).

Modelling All Alternatives also showed some promise in spatial metrics on the small problems for which it can successfully compute objective sets. As such, it could be used in models with very few MGA dimensions. However, MAA's dimensional limitations are not trivial, and do meaningfully restrict its utility.

Generally, HSJ MGA should not be used for capacity expansion models. For HSJ, its lack of exploratory power makes generating useful and generalizable takeaways from the solutions found difficult if not impossible, and the prospects for increasing its power seem limited due to its lack of parallelizability.

Modelling to Generate Alternatives is a powerful methodology to generate a wide range of near-optimal feasible solutions for energy system planning problems, with many possible applications to enhance decision support. This paper has shown, however, that the vector selection methodology used to explore the near-optimal feasible space has profound implications for the utility of any MGA application. While geometrically calculated methods like Modelling All Alternatives provide reasonable spatial exploration in tests, it is restricted by computational difficulties to relatively small-scale problems exploring variation in no more than a couple dozen decision variables (or groupings of variables), limiting its practical applicability in macro-energy system modelling, where higher dimensional MGA is often desirable depending on the research question. Hop-Skip-Jump MGA, on the other hand, performs well in a computational run-time sense, but has extremely limited exploratory power, limiting the insights it can deliver, and potentially giving a false sense of the tradeoffs available within a given circumstance and cost slack. The best traditional methods, per the testing conducted here, are the heuristic-type Random Vector and Capacity Min/Max methods, which performed well in both volumetric expansion and runtime in both a contrived testbed setting and in actual macro-energy systems modelling. However, each of these methods presents the opposite side of a tradeoff: Random Vector offers great guarantees that all dimensions will be explored, whether they are included in the MGA objective function explicitly or not, but struggles to find each variable's most extreme values, while Capacity Min/Max ensures that the dimensions that are specifically included in the MGA objective function are explored close to their extremes but may not optimize every dimension. Notably, as demonstrated here, these methods can be improved via their combination, capturing the strengths of each while mitigating their relative weaknesses. For future applications of MGA, it is thus recommended that researchers use a multi-threaded combination of the Random Vector and Capacity Min/Max methods (or similar strategies) to improve decision support and explore the near-optimal feasible space of their macro-energy system models.